\documentclass[12pt]{article}
\usepackage{amsmath}
\usepackage{amssymb}

\usepackage[cp1251]{inputenc}
\usepackage[russian]{babel}
\newcommand{\il}[2]{\int\limits_{#1}^{#2}}

\newcommand{\ph}{\phantom{a}}
\newcommand{\phh}{\phantom{aaa}}

\newcommand{\sist}[2]{\left\{
\begin{array}{l}
{#1}\\
\ph\\
{#2}
\end{array}
\right.}

\begin{document}
MSC  34L30,   34C99

\vskip 20pt
{\bf \centerline {A new global solvability criterion }

\centerline{for matrix Riccati equations}

\vskip 20pt \centerline {G. A. Grigorian}}
\centerline{\it Institute  of Mathematics NAS of Armenia}
\centerline{\it E -mail: mathphys2@instmath.sci.am}
\vskip 20 pt

Abstract. We use a new approach with a matrix transformation to obtain a new global solvability criterion for matrix Riccati equations. The proven theorem completes an well known result in directions of extension of classes of coefficient of equations and of extension of initial values,under which the solutions of equations are continuable to $\infty$.

\vskip 20pt
Key words: matrix Riccati equations,  nonnegative (positive) definiteness, Hermitian matrices, The Liouville's formula.

\vskip 30pt

{\bf 1. Introduction}. Let $P(t), \ph Q(t), \ph R(t)$ and $S(t)$ be complex-valued locally integrable  matrix functions of dimension $n\times n$ on $[t_0,\infty)$. Consider the matrix Riccati equation
$$
Z' + Z P(t) Z + Q(t) Z + Z R(t) + S(t) = 0, \phh t \ge t_0. \eqno (1.1)
$$
It follows from the general theory of normal systems of ordinary differential equations for any matrix $M$ of dimension $n\times n$ and $t_1 \ge t_0$ there exists $t_2 > t_1 \ph (t_2 \le \infty)$ such that the solution $Z(t)$ of the last equation with $Z(t_1) = M$ exists on $[t_1,t_2)$ and is the unique. The main interest from the point of view of qualitative theory of differential equations represents  the case $t_2 = \infty$.

\vskip 10pt

{\bf Remark 1.1.} {\it Criteria for existence of regular solutions of Eq. (1.1) are obtained in [1-3] (see [1], Theorems 3.5, 3.6, [2,3]).}

\vskip 10pt

For any Hermitian matrices $H_1$ and $H_2$ (of the same dimension) we denote by $H_1 \ge H_2 \ph (>0)$ the nonnegative (positive) definiteness of $H_1 - H_2$, by the symbol $H_1 \le 0 \ph (<0)$ we denote the nonnegative (positive) definiteness of $-H_1$. The next theorem describes a wide class of Eq. (1.1), having solutions on $[t_0,\infty)$.

\vskip 10pt

{\bf Theorem 1.1. ([1, Theorem 3.5])}. {\it If $P(t) \ge 0, \ph S(t) \le 0, \ph R(t) = Q^*(t), \ph t \ge t_0$, then the unique solution $Z(t)$ of Eq. (1.1) with piecewise continuous and locally integrable  coefficients exists for $t \ge t_0$ with
$$
0\le Z(t) \le \widetilde{Z}(t), \phh t \ge t_0.
$$
where $\widetilde{Z}(t)$ is the solution of the linear equation
$$
Z' + A^*(t) Z + Z A(t) + S(t) = 0
$$
with $\widetilde{Z}(t_0) = Z_0 \ge 0$.}

\vskip 10pt

Below we obtain new conditions for the  existence of global  solutions of Eq. (1.1).

\vskip 10pt

{\bf 2. Auxiliary propositions}. The next three  lemmas have importance  in the proof of the existence of global  solutions od Eq. (1.1). Let $M_l \equiv (m_{ij}^l)_{i,j=1}^n, \ph l=1,2$ be complex-valued matrices.

\vskip 10pt

{\bf Lemma 2.1.} {\it The equality
$$
tr(M_1 M_2) = tr (M_2 M_1)
$$
is valid.}

Proof.  We have $tr (M_1 M_2) = \sum\limits_{j=1}^n(\sum\limits_{k=1}^n m_{jk}^1 m_{kj}^2) = \sum\limits_{k=1}^n(\sum\limits_{j=1}^n m_{jk}^1 m_{kj}^2) = \sum\limits_{k=1}^n(\sum\limits_{j=1}^n m_{kj}^2 m_{jk}^1) = tr (M_2 M_1).$ The lemma is proved.

\vskip 10pt

{\bf Lemma 2.2.} {\it Let $H_j, \ph j=1,2$ be Hermitian matrices such that $H_j\ge 0, \ph j=1,2.$ Then
$$
tr (H_1 H_2) \ge 0.
$$
}

Proof. Let $U$ be an unitary transformation such that $\widetilde{H}_1 \equiv U H_1 U^* = diag \{h_1,\dots,h_n\}$.
Since any unitary transformation preserves the nonnegative definiteness of any Hermitian  matrix    we have
$$
h_j \ge 0, \phh j=\overline{1,n}. \eqno (2.1)
$$
Let $\widetilde{H}_2 \equiv U H_2 U^* = (h_{ij})_{i,j=1}^n.$ As for as $\widetilde{H}_2$ is Hermitian it follows that (see [4], p. 300, Theorem 20)  $h_{jj} \ge 0, \ph j=\overline{1,n}$. This together with (2.1) implies
$$
tr (H_1 H_2) = tr([U H_1 U^*] [U H_2 U^*]) = tr (\widetilde{H}_1 \widetilde{H}_2) = \sum\limits_{j=1}^n h_j h_{jj} \ge 0.
$$
The Lemma is proved.

\vskip 10pt

{\bf Lemma 2.3.} {\it Let $H \ge 0$ be a Hermitian matrix of dimension $n \times n$ and let $V$ be any matrix of the same dimension. Then
$$
V H V^* \ge 0. \eqno (2.2)
$$
}

Proof. For any vectors $x$ and $y$ of dimension $n$ denote by $\langle x,y\rangle$ their scalar product. Then $\langle V H V^* x, x\rangle = \langle H(V^* x), (V^* x)\rangle \ge 0$ (since $H \ge 0$). Hence (2.2) is valid. The lemma is proved.

Along with Eq. (1.1) consider the linear matrix system
$$
\sist{\Phi' = R(t) \Phi + P(t) \Psi,}{\Psi' = -S(t)\Phi - Q(t) \Psi,} \ph t \ge t_0. \eqno (2.3)
$$
By a solution of this system we mean a pair $(\Phi(t),\Psi(t))$ of absolutely continuous matrix functions of dimension $n\times n$ on $[t_0,\infty)$, satisfying (2.3) almost everywhere on $[t_0,\infty)$. It is not difficult to verify that all solutions of Eq. (1.1), existing on $[t_1,t_2) \subset [t_0,\infty)$ are connected with solutions $(\Phi(t),\Psi(t))$ of the system (2.3) by relations
$$
\Phi'(t) = [R(t) + P(t)Z(t)]\Phi(t), \phh \Psi(t) = Z(t) \Phi(t), \phh t \in [t_1,t_2).  \eqno (2.4)
$$
By the Liouville's formula from here it follows
$$
\det \Phi(t) = \det \Phi(t_0)\exp\biggl\{\il{t_1}{t} tr \Bigl[R(\tau) + P(\tau) Z(\tau)\Bigr] d \tau\biggr\}, \phh t\in [t_1,t_2), \eqno (2.5)
$$
$$
\overline{\det \Phi(t)} = \overline{\det \Phi(t_0)}\exp\biggl\{\il{t_1}{t} tr \Bigl[R^*(\tau) + Z^*(\tau)P(\tau) \Bigr]d\tau\biggr\}, \phh t\in [t_1,t_2). \eqno (2.6)
$$

\vskip 10pt

{\bf 3. The main result}.

\vskip 5pt
{\bf Definition 3.1.} {\it An interval $[t_1,t_2) \subset [t_0,\infty)$ is called the maximum existence interval for a solution $Z(t)$ of Eq. (1.1), if $Z(t)$ exists on $[t_1,t_2)$ and cannot be continued to the right from $t_2$ as a solution of Eq. (1.1).
}

Let $U(t)$ be a complex-valued absolutely continuous $n\times n$ matrix function on $[t_0,\infty)$ such that $\det U(t) \ne 0.$
For any absolutely continuous matrix function $\Lambda(t)$ of dimension $n\times n$ on $[t_0,\infty)$  we set
$$
S_{U,\Lambda}(t)\equiv \Lambda'(t) + \Lambda(t)  P(t) U(t) \Lambda(t) + \Bigl[U^{-1}(t) U'(t)  + U^{-1}(t)Q(t)U(t)\Bigr] \Lambda(t) + \phantom{aaaaaaaaaaaaaaaaaaaa}
$$
$$
 \phantom{aaaaaaaaaaaaaaaaaaaaaaaaaaaaaaaaaaaaaaaaaaaaaaa}+ \Lambda(t) R(t) + U^{-1}(t)S(t), \ph t \ge t_0.
$$
 \vskip 10pt

{\bf Theorem 3.1.} {\it Let the following  conditions be satisfied:

\noindent
for some absolutely continuous matrix function $\Lambda(t)$ of dimension $n\times n$ on $[t_0,\infty)$

\noindent
I)  $P(t) U(t)$ is Hermitian and     $P(t) U(t) \ge 0, \ph t \ge t_0,$

\noindent
II)   $R(t) = (U^{-1}(t)Q(t)U(t))^* + P(t) U(t)[\Lambda^*(t) - \Lambda(t)]  +  (U^{-1}(t) U'(t))^* +\mu(t) I , \ph t \ge t_0$, where $\mu(t)$ is a real-valued locally integrable function on $[t_0,\infty)$,

\noindent
III) $S_{U,\Lambda} (t) + S_{U,\Lambda} ^*(t) \le 0, \ph t\ge t_0.$

\noindent
Then every solution $Z(t)$ of Eq. (1.1) with $U^{-1}(t_0)Z(t_0) + Z^*(t_0)(U^*)^{-1}(t_0) > \Lambda(t_0) + \Lambda^*(t_0)$  exists on $[t_0,\infty)$   and
$$
U^{-1}(t)Z(t) + Z^*(t)(U^*)^{-1}(t) > \Lambda(t) + \Lambda^*(t_0), \phh t \ge t_0.  \eqno (3.1)
$$
Furthermore, if $U^{-1}(t_0) + (U^{-1}(t_0))^*> 0$, then every solution $Z(t)$ of Eq. (1.1) with $U^{-1}(t_0)Z(t_0) + Z^*(t_0)(U^*)^{-1}(t_0) \ge \Lambda(t_0) + \Lambda^*(t_0)$  exists on $[t_0,\infty)$   and
$$
U^{-1}(t)Z(t) + Z^*(t)(U^*)^{-1}(t) \ge \Lambda(t) + \Lambda^*(t_0), \phh t \ge t_0.  \eqno (3.2)
$$

}

Proof.  In Eq. (1.1) we substitute
$$
Z = U(t)W, \phh t \ge t_0. \eqno (3.3)
$$
We obtain
$$
W' + WP(t)U(t)W + \Bigl(U^{-1}(t) U'(t) + U^{-1}(t)Q(t)U(t)\Bigr)W + W R(t) + U^{-1}(t) S(t) = 0, \ph t \ge t_0.
$$
The next substitution
$$
W = L + \Lambda(t) \eqno (3.4)
$$
reduces the last equation to the following
$$
L' + LP(t)U(t)L + Q_{U,\Lambda}(t)L + L R_{U,\Lambda}(t) + S_{U,\Lambda}(t) = 0, \phh t \ge t_0.  \eqno (3.5)
$$
where
$$
Q_{U,\Lambda}(t)\equiv U^{-1}(t) U(t) + U^{-1}(t)Q(t)U(t) +\Lambda(t)P(t)U(t), \ph R_{U,\Lambda}(t)\equiv R(t) +  P(t)U(t)\Lambda(t),
$$
$t \ge t_0.$
Let $Z(t)$ be a solution of Eq. (1.1) with $U^{-1}(t_0)Z(t_0) + Z^*(t_0)(U^{-1}(t_0))^* > \Lambda(t_0) + \Lambda^*(t_0)$, and let $[t_0,t_1)$ be its maximum existence interval. Then by (3.3) and (3.4) $L(t) \equiv U^{-1}(t) Z(t) - \Lambda(t), \ph t \in [t_0,t_1)$ is a solution of Eq. (3.5) on $[t_0,t_1)$ and  $[t_0,t_1)$  is its maximum existence interval. Obviously,  $L(t_0) + L^*(t_0) > 0$.

Show that
$$
L(t) + L^*(t) > 0, \phh t \in[t_0,t_1). \eqno (3.6)
$$
Suppose this is not so. Then there exists $t_2 \in (t_0,t_1)$ such that
$$
L(t) + L^*(t) > 0, \ph t \in [t_0,t_2), \eqno (3.7)
$$
and
$$
\det [L(t_2) + L^*(t_2)] = 0. \eqno (3.8)
$$
Since $L(t)$ is a solution of Eq. (3.5) on $[t_0,t_1)$ we have
$$
L'(t) + L(t) P(t)U(t) L(t) + Q_{U,\Lambda}(t) L(t) + L(t) R_{U,\Lambda}(t) + S_{u,\Lambda}(t) = 0,
$$
$$
[L^*(t)]' + L^*(t)P(t)U(t) L^*(t) + R_{u,\Lambda}^*(t) Z^*(t) + Z^*(t) Q_{U,\Lambda}^*(t) + S_{U,\Lambda}^*(t) = 0,
$$
almost everywhere on $[t_1,t_2)$. Summing up these equalities, taking into account that $P(t)U(t), \ph t \ge t_0$ is Hermitian and making some simplifications  we obtain
$$
[L(t) + L^*(t)]' + [L(t) + L^*(t)] P(t)U(t) [L(t) + L^*(t)] + \frac{[Q_{U\Lambda}(t) + R_{U\Lambda}^*(t)]}{2}[L(t) + L^*(t)] +
$$
$$
+[L(t) + L^*(t)] \frac{[Q_{U,\Lambda}^*(t) + R_{U,\Lambda}(t)]}{2} +\frac{[Q_{U,\Lambda}(t) - R_{U,\Lambda}^*(t)]}{2}[L(t) - L^*(t)] +
$$
$$
+[L(t) - L^*(t)] \frac{[R_{U,\Lambda}(t) - Q_{U,\Lambda}^*(t)]}{2} +
$$
$$
+[S_{U,\Lambda}(t) + S_{u,\Lambda}^*(t)] - L(t) P(t)U(t) L^*(t) - L^*(t)  P(t)U(t) L(t) = 0 \eqno (3.9)
$$
almost everywhere on $[t_0,t_1)$. It follows from (3.7) that $(L(t) + L^*(t))^{-1}$ exists on $[t_0,t_2)$. Then (3.9) allows  to interpret $ L(t) + L^*(t), \ph t\in [t_0,t_2)$ as a solution of the following linear matrix differential equation
$$
\mathcal{L}' + \biggl\{[L(t) + L^*(t)]  P(t)U(t) + \frac{Q_{U,\Lambda}(t) + R_{U,\Lambda}^*(t)}{2} + [L(t) + \phantom{aaaaaaaaaaaaaaaaaaaaaaaaaaaaaaaaaaaaaaaaaa}
$$
$$
+L^*(t)]\frac{Q_{U,\Lambda}^*(t) + R_{U,\Lambda}(t)}{2}[L(t) + L^*(t)]^{-1} +
$$
$$
+\Bigl(\frac{Q_{U,\Lambda}(t) - R_{U,\Lambda}^*(t)}{2}[L(t) - L^*(t)] + [L(t) - L^*(t)]\frac{R_{U,\Lambda}(t) - Q_{u,\Lambda}^*(t)}{2}\Bigr) [L(t) + L^*(t)]^{-1}+
$$
$$
+\Bigl(S_{U,\Lambda}(t) + S_{U,\Lambda}^*(t) - L(t) P(t) U(t) L^*(t) - L^*(t) P(t)U(t) L(t)\Bigr) [L(t) + L^*(t)]^{-1}\biggr\} \mathcal{L} = 0,
$$
$t \in [t_0,t_2).$
Then by the Liouville formula we have
$$
\det [L(t) + L^*(t)] = \det [L(t_0) + L^*(t_0)] \exp\biggl\{ -\il{t_0}{t} tr  \biggl[[L(\tau) + L^*(\tau)]  P(\tau) U(t) +
$$
$$
+ \frac{Q_{U,\Lambda}(\tau) + R_{U,\Lambda}^*(\tau)}{2} + [L(\tau) + L^*(\tau)]\frac{Q_{U,\Lambda}^*(\tau) + R_{U,\Lambda}(\tau)}{2}[L(t) + L^*(\tau)]^{-1}+
$$
$$
\Bigl(\frac{Q_{U,\Lambda}(\tau) - R_{U,\Lambda}^*(\tau)}{2}[L(\tau) - L^*(\tau)] + [L(\tau) - L^*(\tau)]\frac{R_{U,\Lambda}(\tau) - Q_{U,\Lambda}^*(\tau)}{2}\Bigr) [L(\tau) + L^*(\tau)]^{-1}-
$$
$$
+\Bigl(S_{U,\Lambda}(\tau) + S_{U,\Lambda}^*(\tau) - \phantom{aaaaaaaaaaaaaaaaaaaaaaaaaaaaaaaaaaaaaaaaaaaaaaaaaaaaa}
$$
$$
-L(\tau) P(\tau)U(t) L^*(\tau) - L^*(\tau) P(\tau)U(t) L(\tau)\Bigr) [L(\tau) + L^*(\tau)]^{-1}\biggr]d\tau\biggr\}, \ph t\in [t_0,t_2). \eqno (3.10)
$$
By the condition II) we have
$$
\frac{Q_{U,\Lambda}(t) - R_{U,\Lambda}^*(t)}{2}(L(t) - L^*(t))   + (L(t) - L^*(t))\frac{R_{U,\Lambda}(t) - Q_{U,\Lambda}^*(t)}{2} = \phantom{aaaaaaaaaaaaaaaaaaaaaaaaaaa}
$$
$$
\phantom{aaaaaaaaaaaaaaaaaaaaaaaaaaaaaaaaaaa}= \frac{\mu(t) - \overline{\mu(t)}}{2}(L(t) - L^*(t)) = 0, \ph t \in [t_0,t_2)
$$
Therefore,
$$
tr \biggl[\biggl(\frac{Q_{U,\Lambda}(t) - R_{U,\Lambda}^*(t)}{2}(L(t) - L^*(t)) + \phantom{aaaaaaaaaaaaaaaaaaaaaaaaaaaaaaaaaaaaaaaaaaaaaaaaaaaaaa}
 $$
 $$
 \phantom{aaa} +(L(t) - L^*(t))\frac{R_{U,\Lambda}(t) - Q_{U,\Lambda}^*(t)}{2}\biggr)(L(t) + L^*(t))^{-1}\biggr] = 0,  \ph t\in[t_0,t_2). \eqno (3.11)
$$
By Lemma 2.3 it follows from the condition I)  that
$$
L(t) P(t) U(t) L^*(t) + L^*(t) P(t)U(t) L(t) \ge 0, \ph t \in [t_0,t_2].
$$
Then since $(L(t) + L^*(t))^{-1} > 0, \ph t \in [t_0,t_2)$ by Lemma 2.2 we have
$$
tr \Bigl[(L(t)P(t)U(t)L^*(t) + L^*(t)P(t)U(t)L(t))(L(t) + L^*(t))^{-1}\Bigr] \ge 0, \ph t \in [t_0,t_2). \eqno (3.12)
$$
By the same reason from the condition III) of the theorem we obtain
$$
tr\Bigl[(S_{U,\Lambda}(t) + S_{U,\Lambda}^*(t))(L(t) + L^*(t))^{-1}\Bigr]\le 0, \phh t \in [t_0,t_2). \eqno (3.13)
$$
Since $L(t) + L^*(t), \ph t \in [t_0,t_1)$ is absolutely continuous on $[t_0,t_2]$ we have
$$
\biggl|\il{t_0}{t_2}tr \Bigl[(L(t) + L^*(t))u(t)P(t)\Bigr]d t\biggr| \le c =const. \eqno (3.14)
$$
Obviously
$$
tr\Bigl[(L(t) + L^*(t))\frac{Q_{U,\Lambda}^*(t) + R_{U,\Lambda}(t)}{2}(L(t) + L^*(t))^{-1}\Bigr] = tr \biggl[\frac{Q_{U,\Lambda}^*(t) + R_{U,\Lambda}(t)}{2}\biggr], \ph t\in [t_0,t_2).
$$
This together with  (3.10) - (3.14)   implies that
$$
\det (L(t_2) + L^*(t_2)) \ge \det (L(t_0) + L^*(t_0)) \exp\bigl\{-c\} > 0,
$$
which contradicts (3.8). The obtained contradiction proves (3.6). By (3.3) and (3.4)
it follows from (3.6) that $Z(t)$ exists on $[t_0,t_1)$ and
$$
U^{-1}(t) Z(t) + Z^*(t)(U^{-1}(t))^* > \Lambda(t) + \Lambda^*(t), \phh t \in [t_0,t_1). \eqno (3.15)
$$
Let $(\Phi(t),\Psi(t))$ be a solution of the system (3.4)  with $\Phi(t_0) = I, \ph \Psi(t_0) = Z(t_0)$. Then by   (2.5), (2.6)  and Lemma 2.1  we have
$$
|\det\Phi(t)|^2 = |\det\Phi(t_0)|^2 \exp\biggl\{\il{t_0}{t} tr \biggl[R(\tau) + R^*(\tau) + P(\tau)(Z(\tau) + Z^*(\tau))\biggr]d\tau\biggr\}, \ph t \in [t_0,t_1).
$$
By Lemma 2.2 from the condition $I)$ of the theorem we derive $\il{t_0}{t}te \Bigl[P(\tau)(Z(\tau) + Z^*(\tau))\Bigr] d \tau = \linebreak \il{t_0}{t}te \Bigl[u(\tau)P(\tau)\bigl[Z(\tau)/u(\tau) + Z^*(\tau)/\overline{u}(\tau)\bigr]\Bigr] d \tau \ge~ 0, \ph t \in [t_0,t_1)$. Then from the last equality we obtain.
$$
|\det\Phi(t_1)|^2 \ge |\det\Phi(t_0)|^2 \exp\biggl\{\min\limits_{t\in[t_0,t_1]}\il{t_0}{t} tr \biggl[R(\tau) + R^*(\tau) \biggr]d\tau\biggr\}> 0.
$$
Hence, $\det \Phi(t) \ne 0, \ph t\in [t_0, t_1 +\varepsilon)$ for some $\varepsilon > 0$. By (3.5) from here it follows that $Z_1(t)\equiv \Psi(t)\Phi^{-1}(t), \ph t \in[t_0,t_1+\varepsilon)$ is a solution of Eq. (1.1) on $[t_0,t_1 +\varepsilon)$ and coincides with $Z(t)$ on $[t_0,t_1)$. By the uniqueness theorem it follows from here that $[t_0,t_1)$ is not the maximum existence interval for $Z(t)$. Therefore $t_1 = \infty$. This together with (3.15) implies (3.1).
 Let us prove  (3.2).
For any $\varepsilon > 0$ denote by $Z_\varepsilon(t)$ the solution of Eq. (1.1), satisfying the initial condition
$$
Z_\varepsilon(t_0) = Z(t_0) +\varepsilon I,
$$
where $Z(t)$ is a solution of Eq. (1.1), with  $U^{-1}(t_0) Z(t_0) + Z^*(t_0)(U^{-1}(t_0))^* \ge  \Lambda(t_0) + \Lambda^*(t_0)$. Obviously
$$
U^{-1}(t_0) Z_\varepsilon(t_0) + Z_\varepsilon ^*(t_0)(U^{-1}(t_0))^* \ge  \varepsilon [U^{-1}(t_0) + (U^{-1}(t_0))^*] + \Lambda(t_0) + \Lambda^*(t_0) >  \Lambda(t_0) + \Lambda^*(t_0).
$$
Then by already proven $Z_\varepsilon(t)$ exists on $[t_0,\infty)$ and
$$
U^{-1}(t) Z_\varepsilon(t) + Z_\varepsilon^*(t)(U^{-1}(t))^* > \Lambda(t_0) + \Lambda^*(t_0), \phh t \ge t_0. \eqno (3.16)
$$
Let $[t_0,T)$ be the maximum existence interval for $Z(t)$. Show that
$$
T = \infty. \eqno (3.17)
$$
Suppose $T< \infty$. Let $\lambda(t)$ be the least eigenvalue of $U^{-1}(t)Z(t) + Z^*(t)(U({-1}(t))^* - (\Lambda(t) + \Lambda^*(t))$ and $\lambda_\varepsilon(t)$ be the least eigenvalue of $U^{-1}(t)Z_\varepsilon(t) + Z_\varepsilon^*(t)(U^{-1}(t))^* - (\Lambda(t) + \Lambda^*(t)), \ph t \in [t_0,T)$. It follows from (3.16) that
$$
\lambda_\varepsilon(t) > 0, \phh t \ge t_0. \eqno (3.18)
$$
Since the solutions of Eq. (1.1) are continuously dependent on their initial values we have $\lambda_\varepsilon(t) \to \lambda (t)$ as $\varepsilon \to 0$. This together with (3.18)  implies that $\lambda(t) \ge 0,\ph t \in [t_0,T).$ Therefore
$$
U^{-1}(t)Z(t) + Z^*(t)(U^{-1}(t))^* \ge  \Lambda(t) + \Lambda^*(t), \phh t\in [t_0,T). \eqno (3.19)
$$
Let $(\Phi_0(t),\Psi_0(t))$ be a solution of the system (2.3) with $\Phi_0(t_0) = I, \ph \Psi_0(t_0) = Z(t_0)$. Then by (2.5) and (2.6) we have
$$
|\det\Phi_0(t)|^2 = |\det\Phi_0(t)|^2 \exp\biggl\{\il{t_0}{t} tr \biggl[R(\tau) + R^*(\tau) + P(\tau)(Z(\tau) + Z^*(\tau))\biggr]d\tau\biggr\}=
$$
$$
=|\det\Phi_0(t)|^2 \exp\biggl\{\il{t_0}{t} tr \biggl[R(\tau) + R^*(\tau) + u(\tau)P(\tau)(Z(\tau)/u(\tau) + Z^*(\tau)/\overline{u}(\tau))\biggr]d\tau\biggr\},
$$
$t \in [t_0,T).$
By Lemma 2.2 it follows from here and (3.19)
that
$$
|\det \Phi_0(T)|^2 \ge |\det \Phi_0(t_0)|^2\exp\biggl\{\min\limits_{t \in [t_0,T]}\il{t_0}{t} tr \Bigl[R(\tau) + R^*(\tau)\Bigr]d\tau\biggr\} > 0.
$$
Hence, $\det \phi_0(t) \ne 0, \ph t \in [t_0,T+\delta)$ for some $\delta > 0$. Then by (2.6)  $Z_0(t) \equiv \Psi_0(t)\Phi_0^{-1}(t), \ph t \in [t_0,T+\delta)$ is a solution of Eq. (1.1) on $[t_0,T+\delta)$, which coincides with $Z(t)$ on $[t_0,T)$. By the uniqueness theorem it follows from here that $[t_0,T)$ is not the maximum existence interval for $Z(t)$, which contradicts our supposition. The obtained contradiction proves (3.17). From (3.17) and (3.19) it follows (3.2). The theorem is proved.

\vskip 10pt

{\bf Remark 5.1.} {\it Theorem 1.1 is a complement  of Theorem 3.1 in two directions:

\noindent
I)  the set of coefficients of Eq. (1.1) is extended;

\noindent
II) the set of initial values, for which the solutions of Eq. (1.1) are continuable to $\infty$, is extended.
}

\vskip 10 pt

If $P(t)$ is absolutely continuous and $\det P(t) \ne 0, \ph t \ge t_0$, then, obviously, $P(t)P^*(t)>~ 0, \linebreak t \ge t_0$ and is Hermitian. This together with Theorem 3.1 implies immediately

\vskip 10pt
{\bf Corollary 3.1.} {\it Let $P(t)$ be absolutely continuous on $[t_0,\infty)$ such that $\det P(t) \ne 0  \ph t \ge t_0$. If for some absolutely continuous matrix function $\Lambda(t)$ of dimension $n\times n$ on $[t_0,\infty)$ and a  real-valued locally integrable function $\mu(t)$ on $[t_0,\infty)$

\noindent
$R(t) = (P^*(t)Q^*(t)(P^{-1}(t))^* + P(t)P^*(t)[\Lambda^*(t) - \Lambda(t)]  +   (P'(t))^*P^{-1}(t) +\mu(t) I , \ph t \ge t_0$ and

\noindent
 $S_{P^*,\Lambda} (t) + S_{P^*,\Lambda} ^*(t) \le 0, \ph t\ge t_0.$

\noindent
Then every solution $Z(t)$ of Eq. (1.1) with $(P^{-1}(t_0))^*Z(t_0) + Z^*(t_0)P^{-1}(t_0) > \Lambda(t_0) + \Lambda^*(t_0)$  exists on $[t_0,\infty)$   and
$$
(P^{-1}(t))^*Z(t) + Z^*(t)P^{-1}(t) > \Lambda(t) + \Lambda^*(t)
 \phh t \ge t_0.
$$
Furthermore, if $P^{-1}(t_0) + (P^{-1}(t_0))^*> 0$, then every solution $Z(t)$ of Eq. (1.1) with   $(P^{-1}(t_0))^*Z(t_0) + Z^*(t_0)P^{-1}(t_0) \ge \Lambda(t_0) + \Lambda^*(t_0)$
$$
(P^{-1}(t))^*Z(t) + Z^*(t)P^{-1}(t) \ge \Lambda(t) + \Lambda^*(t)
 \phh t \ge t_0.
$$

}

\vskip 20pt

\centerline{\bf References}

\vskip 20pt

\noindent
1. G. Freiling, A survey  of nonsymmetric Riccati equations. Linear Algebra and its \linebreak\phantom{aa}  Applications.  351-352 (2002) pp. 243-270.

\noindent
2. G. Freiling, G. Jank, A. Sarychev, Non-blow-up conditions for Riccati-type matrix \linebreak\phantom{aa} differential and difference equations. Results math. 37 (1998) pp. 84-103.

\noindent
3. H. W. Knobloch, M.Pohl, On Riccati matrix differential equations. Results Math. \linebreak\phantom{aa} 31 (19970, PP. 337-364.

\noindent
4. F. G.  Gantmacher,  Theory of matrices. Moskow, Nauka, 1966, 576 pages.

\end{document}